\documentclass[12pt]{article}
\usepackage{a4wide, amstext, amsmath, amssymb, graphicx, array, epsfig, psfrag}
\usepackage{amstext, amsmath, amssymb, graphicx,color}

\newcommand{\bfJ}{{\boldsymbol J}}
\newcommand{\bfn}{\boldsymbol n}
\newcommand{\bfkappa}{\boldsymbol\kappa}
\newcommand{\bfF}{\boldsymbol F}
\newcommand{\bff}{\boldsymbol f}
\newcommand{\bfp}{\boldsymbol p}
\newcommand{\bfu}{\boldsymbol u}
\newcommand{\bfeps}{\boldsymbol\varepsilon}
\newcommand{\bfsig}{\boldsymbol\sigma}
\newcommand{\bfv}{\boldsymbol v}

\newcommand{\bfI}{\boldsymbol I}
\newcommand{\bfP}{\boldsymbol P}
\newcommand{\bfx}{\boldsymbol x}
\newcommand{\bfq}{\boldsymbol q}

\numberwithin{equation}{section}

\newtheorem{rem}{Remark}[section]

\date{}
\begin{document}
\title{\bf Intrinsic finite element modeling of a linear membrane shell problem}
\author{Peter Hansbo
\footnote{Department of Mechanical Engineering, J\"onk\"oping University, 
SE-55111 J\"onk\"oping, Sweden.} \mbox{ }
Mats G.\ Larson
\footnote{Department of Mathematics and Mathematical Statistics, Ume{\aa} University, SE-90187 Ume{\aa}, Sweden} }
\numberwithin{equation}{section} \maketitle
\begin{abstract}
A Galerkin finite element method for 
the membrane elasticity problem on a meshed surface is constructed by using two-dimensional
elements extended into three dimensions. 
The membrane finite element model is established using the intrinsic approach suggested by Delfour and Zol\'esio \cite{DeZo95}.
\end{abstract}

\section{Introduction}

Models of thin-shell structures are often established using differential geometry
to define the governing differential equations in two dimensions, cf. Ciarlet \cite{Ci00} 
for an overview. A simpler approach is the classical engineering trick of 
viewing the shell as an assembly of flat elements, in which simple transformations of the two-dimensional stiffness matrices are performed, cf., e.g., Zienkiewciz \cite{zienk}.
In contrast to these approaches, Delfour and Zolesio \cite{DeZo95,DeZo96,DeZo97} 
established elasticity models on surfaces using the signed distance function,
which can be used to describe the geometric properties of a surface. In particular, the 
intrinsic tangential derivatives were used for modeling purposes as the main 
differential geometric tool and the partial differential equations were established in three
dimensions. A similar concept had been used earlier in a finite element setting for the 
numerical discretization of the Laplace-Beltrami operator on surfaces by Dziuk \cite{Dz88}, 
resulting in a remarkably clean and simple implementation. For diffusion-like problems, 
the intrinsic approach has become the focal point of resent research on numerical solutions 
of problems posed on surfaces, cf., e.g., \cite{BaMoNo04,BaGaNa08,DeDzElHe10,DeDz07,ElSt10,OlReGr09}

The purpose of this paper is to begin to explore the possibilities of the intrinsic approach 
in finite element modeling of thin-shell structures, focusing on the simplest model, that 
of the membrane shell without bending stiffness. We derive a membrane model using the intrinsic framework and generalize the finite element approach of \cite{Dz88}. Finally, we give some 
elementary numerical examples.

\section{The membrane shell model problem}

\subsection{Basic notation}

We begin by recalling the fundamentals of the approach of Delfour and Zolesio \cite{DeZo95,DeZo96,DeZo97}. Let $\Sigma$ be a smooth two-dimensional surface imbedded in ${\mathbb{R}}^3$, with outward pointing normal $\bfn$. 
If we denote the signed distance function relative to $\Sigma$ by $d(\bfx)$, for $\bfx\in \Bbb{R}^3$, fulfilling $\nabla d = \bfn$, we can define the domain occupied by the membrane by
\[
\Omega_t = \{\bfx\in \Bbb{R}^3: \vert d(\bfx) \vert < t/2\},
\]
where $t$ is the thickness of the membrane. The closest point projection $\bfp:\Omega_t \rightarrow \Sigma$
is given by
\[
\bfp(\bfx) = \bfx -d(\bfx)\bfn(\bfx) ,
\]
the Jacobian matrix of which is 
\[
\nabla \bfp = \bfI -d\nabla\otimes\bfn -\bfn\otimes\bfn
\]
where $\bfI$ is the identity and $\otimes$ denotes exterior product.
The corresponding linear projector
$\bfP_\Sigma = \bfP_\Sigma(\bfx)$, onto the tangent plane of $\Sigma$ 
at $\bfx\in\Sigma$, is given by
\[
\bfP_\Sigma := \bfI -\bfn\otimes\bfn ,
\]
and we can then define the surface gradient
$\nabla_\Sigma$ as
\begin{equation}
\nabla_\Sigma := \bfP_\Sigma \nabla .
\end{equation}
The surface gradient thus has three
components, which we shall denote by
\[
\nabla_\Sigma =: \left[\begin{array}{>{\displaystyle}c}
\frac{\partial}{\partial x^\Sigma_1} \\[3mm]
\frac{\partial}{\partial x^\Sigma_2} \\[3mm]
\frac{\partial}{\partial x^\Sigma_3} \
\end{array}\right] .
\]
For a vector valued function $\bfv(\bfx)$, we define the
tangential Jacobian matrix as
\[
\nabla_\Sigma\otimes\bfv :=\left[\begin{array}{>{\displaystyle}c>{\displaystyle}c>{\displaystyle}c}
\frac{\partial v_1}{\partial x^\Sigma_1} &\frac{\partial v_1}{\partial x^\Sigma_2} & \frac{\partial v_1}{\partial x^\Sigma_3} \\[3mm]
\frac{\partial v_2}{\partial x^\Sigma_1} &\frac{\partial v_2}{\partial x^\Sigma_2} & \frac{\partial v_2}{\partial x^\Sigma_3} \\[3mm]
\frac{\partial v_3}{\partial x^\Sigma_1} &\frac{\partial v_3}{\partial x^\Sigma_2} & \frac{\partial v_3}{\partial x^\Sigma_3}
\end{array}\right] 
\]
and the surface divergence $\nabla_{\Sigma}\cdot\bfv := \text{tr}\nabla_\Sigma\otimes\bfv$.
%
%

\subsection{The surface strain and stress tensors}

We next define a surface strain tensor 
\[
\bfeps_{\Sigma}(\bfu) := \frac12\left(\nabla_\Sigma\otimes \bfu + (\nabla_\Sigma\otimes\bfu)^{\rm T}\right),
\]
which is extensively used in \cite{DeZo95,DeZo96,DeZo97}, where it
is employed to derive models of shells based on purely mathematical arguments.

From a mechanical point of view, the problem of using $\bfeps_{\Sigma}(\bfu)$ as a fundamental measure of strain on a surface lies in it not being an in-plane tensor, in that $\bfeps_{\Sigma}(\bfu)\cdot\bfn \neq\bf 0$. The shear strains associated with the out-of-plane direction are typically neglected in mechanical models, but are present in $\bfeps_{\Sigma}(\bfu )$ (cf. Remark \ref{remarkstrain}). 
To obtain an in-plane strain tensor we need to use the projection twice to define
\[
\bfeps^P_{\Sigma}(\bfu) := \bfP_\Sigma\bfeps(\bfu)\bfP_\Sigma ,
\]
which lacks all out-of-plane strain components. For a shell, where plane stress is assumed,
this strain tensor can still be used, since out-of-plane strains do not contribute to the strain energy.

\begin{rem}\label{remarkstrain}
It is instructive to work out the details at a surface point whose surrounding
is tangential to the $x_1x_2$--plane.
In this case $\bfn = (0,0,1)$,
\[
\bfP_\Sigma = \left[\begin{array}{ccc} 1 & 0 & 0 \\ 0 & 1 & 0 \\ 0 & 0 & 0\end{array}\right] ,\quad
\nabla_\Sigma\otimes \bfu = \left[\begin{array}{>{\displaystyle}c>{\displaystyle}c>{\displaystyle}c} 
\frac{\partial u_1}{\partial x_1} & \frac{\partial u_2}{\partial x_1} & \frac{\partial u_3}{\partial x_1} \\[3mm] 
\frac{\partial u_1}{\partial x_2} & \frac{\partial u_2}{\partial x_2} &\frac{\partial u_3}{\partial x_2} \\ 0 & 0 & 0\end{array}\right] ,
\]
\[
\bfeps_\Sigma(\bfu) = \left[\begin{array}{>{\displaystyle}c>{\displaystyle}c>{\displaystyle}c} 
\frac{\partial u_1}{\partial x_1} & \frac12 \left(\frac{\partial u_1}{\partial x_2}+\frac{\partial u_2}{\partial x_1}\right) & \frac12\frac{\partial u_3}{\partial x_1} \\[3mm] 
\frac12 \left(\frac{\partial u_1}{\partial x_2}+\frac{\partial u_2}{\partial x_1}\right) & \frac{\partial u_2}{\partial x_2} & \frac12\frac{\partial u_3}{\partial x_2} \\ \frac12\frac{\partial u_3}{\partial x_1} & \frac12\frac{\partial u_3}{\partial x_2} & 0\end{array}\right],
\]
and
\[
\bfeps^P_\Sigma = \left[\begin{array}{>{\displaystyle}c>{\displaystyle}c>{\displaystyle}c} 
\frac{\partial u_1}{\partial x_1} & \frac12 \left(\frac{\partial u_1}{\partial x_2}+\frac{\partial u_2}{\partial x_1}\right) & 0 \\[3mm] 
\frac12 \left(\frac{\partial u_1}{\partial x_2}+\frac{\partial u_2}{\partial x_1}\right) & \frac{\partial u_2}{\partial x_2} & 0 \\ 0 & 0 & 0\end{array}\right] .
\]
The terms in $\bfeps_\Sigma$ not present in $\bfeps^P_\Sigma$ are shear strains that are typically neglected for thin structures, and it is clear that in our case $\bfeps^P_\Sigma$ is the relevant strain tensor. 
\end{rem}

However, the tensor $\bfeps^P_{\Sigma}$ is rather cumbersome to use directly in a numerical implementation;
it would be much easier to work with $\bfeps_{\Sigma}$ which can
be establish using tangential derivatives. For this reason, we use the fact that there also holds (as is easily confirmed)
\[
\bfeps^P_{\Sigma}(\bfu) = \bfP_\Sigma\bfeps_{\Sigma}(\bfu)\bfP_\Sigma =\frac12\left(\bfP_\Sigma\nabla_\Sigma\otimes\bfu\bfP_\Sigma +(\bfP_\Sigma\nabla_\Sigma\otimes\bfu\bfP_\Sigma)^{\rm T}\right),
\]
and since $\bfn\cdot\bfeps_{\Sigma}(\bfu)\cdot\bfn=0$ we have the following relation:
\[
\bfeps^P_{\Sigma}(\bfu) = \bfeps_{\Sigma}(\bfu) - \left((\bfeps_{\Sigma}(\bfu)\cdot\bfn)\otimes\bfn + \bfn\otimes (\bfeps_{\Sigma}(\bfu)\cdot\bfn)\right) ,
\]
so that, using dyadic double-dot product,
\[
\bfsig:\bfu\otimes\bfv = (\bfsig\cdot\bfu)\cdot\bfv, \quad \bfu\otimes\bfv:\bfsig = \bfu\cdot(\bfv\cdot\bfsig),
\]
where $\bfsig$ is a tensor and $\bfu$, $\bfv$ are vectors, we arrive at
\begin{equation}\label{doubledot}
\bfeps^P_{\Sigma}(\bfu):\bfeps^P_{\Sigma}(\bfv) = \bfeps_{\Sigma}(\bfu):\bfeps_{\Sigma}(\bfv) - 2(\bfeps_{\Sigma}(\bfu)\cdot\bfn)\cdot(\bfeps_{\Sigma}(\bfv)\cdot\bfn),
\end{equation}
which will be used in the finite element implementation below. We also note that there holds 
\begin{equation}\label{trace}
\text{tr}\,\bfeps^P_\Sigma(\bfv) = \nabla_\Sigma\cdot\bfv ,
\end{equation}
where $ \text{tr}\bfeps = \sum_k\varepsilon_{kk}$.

We shall assume an isotropic stress--strain relation,
\[
\bfsig = 2\mu \bfeps + \lambda \text{tr}\bfeps\, \bfI ,
\]
where $\bfsig$ is the stress tensor and $\bfI$ is the identity tensor.
The Lam\'e parameters $\lambda$ and $\mu$ are related to Young's modulus $E$ 
and Poisson's ratio $\nu$ via
\[
\mu = \frac{E}{2(1+\nu)},\quad \lambda =\frac{ E\nu}{(1+\nu)(1-2\nu)} .
\]
For
the in-plane  stress tensor we thus assume
\[
\bfsig_\Sigma^P:=  2\mu \bfeps^P_\Sigma + {\lambda} \text{tr}\bfeps^P_\Sigma\, \bfP_\Sigma ,
\]
in the plane strain case and, in the plane stress case, which is appropriate for a thin membrane,
\begin{equation}\label{constitutive}
\bfsig_\Sigma^P:=  2\mu \bfeps^P_\Sigma + {\lambda_0} \text{tr}\bfeps^P_\Sigma\, \bfP_\Sigma ,
\end{equation}
where
\[
{\lambda_0}:= \frac{2\lambda\mu}{\lambda+2\mu} = \frac{E\nu}{1-\nu^2}.
\]

\subsection{The membrane shell equations}

Consider a potential energy functional given by
\[
\Pi(\bfu_t) := \frac12 \int_{\Omega_t} \bfsig(\bfu_t) :\bfeps(\bfu_t) d\Omega_t -\int_{\Omega_t}\bff_t\cdot\bfu_t
\]
where $\bff_t$ is of the form $\bff_t = \bff\circ \bfp$. 
Under the assumption of small thickness, we have
\[
\int_{\Omega_t} f(\bfx)\, d\Omega_t \approx \int_{-t/2}^{t/2}\int_{\Sigma}f \, d\Sigma dz
\]
and thus
\begin{align*}
\Pi(\bfu_t) \approx \Pi^P_\Sigma (\bfu) {}& := \frac{t}{2} \int_{\Sigma} \bfsig^P_{\Sigma}(\bfu) :\bfeps^P_{\Sigma}(\bfu) d\Sigma -t\int_{\Sigma}\bff\cdot\bfu\, d\Sigma \\
{}&: = \frac{t}{2} (\bfsig^P_{\Sigma}(\bfu),\bfeps^P_{\Sigma}(\bfu))_\Sigma - t (\bff,\bfu)_\Sigma .
\end{align*}

Minimizing the potential energy leads to the variational problem of finding $\bfu \in V$, where $V$ 
is an appropriate Hilbert space which we specify below, such that
\begin{equation}
a_{\Sigma}(\bfu,\bfv) = l_{\Sigma}(\bfv) \quad \forall \bfv \in V
\end{equation}
where, by (\ref{doubledot}) and (\ref{trace}), 
\begin{align}\nonumber
a_{\Sigma}(\bfu,\bfv) &{}  = (2\mu\bfeps^P_\Sigma(\bfu), \bfeps^P_\Sigma(\bfv))_\Sigma+ ({\lambda_0}\, \text{tr}\,\bfeps^P_\Sigma(\bfu),\text{tr}\,\bfeps^P_\Sigma(\bfv))_{\Sigma}\\\nonumber
&{}= (2\mu\bfeps_\Sigma(\bfu), \bfeps_\Sigma(\bfv))_\Sigma-(4\mu\bfeps_\Sigma(\bfu)\cdot\bfn, \bfeps_\Sigma(\bfv)\cdot\bfn)_{\Sigma}+ ({\lambda_0}\, \nabla_\Sigma\cdot\bfu,\nabla_\Sigma\cdot\bfv)_{\Sigma}, 
\end{align}
and $l_{\Sigma}(\bfv) = (\bff,\bfv)_\Sigma$. 
This variational problem formally coincides with the one analyzed in the classical differential geometric setting by Ciarlet and co-workers \cite{CiSa93,CiLo96}, as shown in \cite{DeZo97}.

Splitting the displacement into a normal part $u_n := \bfu\cdot\bfn$ and a tangential part 
$\bfu_t := \bfu-u_n \bfn$ we have the identity 
\[
\bfeps_\Sigma^P(\bfu) = \bfeps_\Sigma^P(\bfu_t)+u_n\nabla\otimes\bfn =  
\bfeps_\Sigma^P(\bfu_t)+ u_n \bfkappa,
 \]
where $\bfkappa = \nabla \otimes \nabla d$ is the Hessian of the distance function $d$, 
cf. \cite{DeZo97}, The bilinear form can therefore also be written in the form
 \begin{align}\nonumber
a_{\Sigma}(\bfu,\bfv) &{}  = (2\mu (\bfeps^P_\Sigma(\bfu_t) + u_n \bfkappa ) , \bfeps^P_\Sigma(\bfv_t) + v_n \bfkappa)_\Sigma
\\
&\qquad + ({\lambda_0} ( \text{tr}\, \bfeps^P_\Sigma(\bfu_t) + u_n \text{tr} \, \bfkappa ) ,\text{tr}\,\bfeps^P_\Sigma(\bfv_t) + v_n \text{tr}\, \bfkappa)_{\Sigma}\\\nonumber
\end{align}
This means that we do not have full ellipticity in our problem. Based on this observation 
we conclude that the natural function space for the variational formulation is 
\[
V = \{\bfv  : v_n \in L_2(\Sigma )\quad\text{and}\quad \bfv_t \in [H^1(\Sigma)]^2\},
\]
cf. \cite{CiLo96}. The loss of ellipticity have consequences for the numerics and we comment on 
this in the numerical examples below.

Since 
\[
 (\bfsig^P_{\Sigma}(\bfu),\bfeps^P_{\Sigma}(\bfu))_\Sigma= (\bfsig^P_{\Sigma}(\bfu),\bfeps_{\Sigma}(\bfu))_\Sigma
 \]
 we find, using Green's formula, the pointwise equilibrium equation
 \begin{equation}\label{equil}
 -\nabla_\Sigma\cdot \bfsig^P_{\Sigma}(\bfu) = \bff\quad\text{in}\; \Sigma ,
 \end{equation}
 which together with the constitutive law (\ref{constitutive}) defines the intrinsic 
 differential equations of linear elasticity on surfaces.


\section{The finite element method}

\subsection{Parametrization}

Let $\mathcal{T}_h:=\{T\}$ be a conforming, shape regular triangulation of {$\Sigma$}, resulting in a discrete surface $\Sigma_h$. 
We shall here consider an
isoparametric parametrization of the surface  (the same idea can however be used for arbitrary parametrizations). In the numerical examples below we use a piecewise linear approximation, meaning that the elements $T$ will be planar.
For the parametrization we wish to define a map $\bfF : (\xi,\eta)\rightarrow (x,y,z)$ from a reference triangle $\hat T$ defined in a local coordinate system $(\xi, \eta)$ to $T$, for all $T$. To this end, we write
$\bfx = \bfx(\xi, \eta)$, where $\bfx=(x,y,z)$
are the physical coordinates on $\Sigma_h$.
For any given parametrization, we can extend it outside the surface by 
defining
\[
\bfx(\xi,\eta,\zeta) = \bfx(\xi, \eta)+\zeta\,\bfn(\xi,\eta)
\]
where $\bfn$ is the normal and $-t/2 \leq \zeta \leq t/2$. In some models, where the surface is an idealized thin structure, it is natural to think of $t$ as a thickness.  

For the representation of the geometry, we first introduce the following approximation 
of the normal: 
\[
\bfn \approx \bfn^h = \frac{\bfn^h_0}{\vert\bfn^h_0\vert}, \quad \bfn^h_0 = \sum_i\bfn_i\varphi_i(\xi,\eta),
\] 
where $\varphi_i(\xi,\eta)$ are the finite element shape functions on the reference element 
(assumed linear in this paper), and $\bfn_i$ denotes the normals in the nodes of the mesh.
We then consider parametrizations of the type
\begin{equation}\label{approx}
\bfx (\xi,\eta,\zeta) \approx \bfx^h(\xi,\eta,\zeta)= \sum_i\left(\bfx_i\varphi_i(\xi,\eta)+\zeta\,\bfn_i\varphi_i(\xi,\eta)\right)
\end{equation}
where $\bfx_i$ are the physical location of the nodes on the surface. 
For the approximation of the solution, we use a constant extension,
\begin{equation}
\bfu \approx\bfu^h = \sum_i\bfu_i\varphi_i(\xi,\eta)
\end{equation}
where $\bfu_i$ are the nodal displacements,
so that the finite element method is, in a sense, superparametric. Note that only the in-plane variation of the approximate solution will matter since we are looking at in-plane stresses and strains. 
We employ the usual finite element approximation of the physical derivatives of the chosen basis $\{\varphi_i\}$ on the surface, at $(\xi,\eta)$, as
\[
\left[\begin{array}{>{\displaystyle}c}
\frac{\partial \varphi_j}{\partial x}\\[2mm]
\frac{\partial \varphi_j}{\partial y}\\[2mm]
\frac{\partial \varphi_j}{\partial z}\end{array}\right] = \bfJ^{-1}(\xi,\eta,0) \left[\begin{array}{>{\displaystyle}c}
\frac{\partial \varphi_j}{\partial \xi}\\[2mm]
\frac{\partial \varphi_j}{\partial \eta}\\[2mm]
\frac{\partial \varphi_j}{\partial \zeta}\end{array}\right]_{\zeta=0} \quad\text{where}
\quad
\bfJ(\xi,\eta,\zeta) := \left[\begin{array}{>{\displaystyle}c>{\displaystyle}c>{\displaystyle}c}
\frac{\partial x^h}{\partial \xi} & \frac{\partial y^h}{\partial \xi} & \frac{\partial z^h}{\partial \xi}\\[2mm]
\frac{\partial x^h}{\partial \eta} & \frac{\partial y^h}{\partial \eta} & \frac{\partial z^h}{\partial \eta}\\[2mm]
\frac{\partial x^h}{\partial \zeta} & \frac{\partial y^h}{\partial \zeta} & \frac{\partial z^h}{\partial \zeta}\end{array}\right] ,
\]
%
This gives, at $\zeta=0$,
\[
\left[\begin{array}{>{\displaystyle}c}
\frac{\partial \varphi_i}{\partial x}\\[2mm]
\frac{\partial \varphi_i}{\partial y}\\[2mm]
\frac{\partial \varphi_i}{\partial z}\end{array}\right] = \bfJ^{-1}(\xi,\eta,0) \left[\begin{array}{>{\displaystyle}c}
\frac{\partial \varphi_i}{\partial \xi}\\[2mm]
\frac{\partial \varphi_i}{\partial \eta}\\[2mm]
0\end{array}\right] .
\]
With the approximate normals we explicitly obtain
\[
\left.\frac{\partial \bfx^h}{\partial \zeta}\right|_{\zeta=0}=\bfn^h ,
\]
so
\[
\bfJ(\xi,\eta,0) := \left[\begin{array}{>{\displaystyle}c>{\displaystyle}c>{\displaystyle}c}
\frac{\partial x^h}{\partial \xi} & \frac{\partial y^h}{\partial \xi} & \frac{\partial z^h}{\partial \xi}\\[3mm]
\frac{\partial x^h}{\partial \eta} & \frac{\partial y^h}{\partial \eta} & \frac{\partial z^h}{\partial \eta}\\[3mm]
n^h_y & n^h_y & n^h_z\end{array}\right] .
\]
\begin{rem}\label{dziuk}
The approach by Dziuk \cite{Dz88} (and also the classical engineering approach, \cite{zienk}) is, in our setting, a constant-by-element extension of 
the geometry using facet triangles $\{T\}$ so that,
with $\bfn_T$ the normal to the facet,
$\left.\bfx^h (\xi,\eta,\zeta)\right|_T =  \sum_{i}\bfx_i\left.\varphi_i(\xi,\eta)\right|_T+\zeta\bfn_T$, and 
\[
\left[\begin{array}{>{\displaystyle}c}
\frac{\partial \varphi_i}{\partial x}\\[3mm]
\frac{\partial \varphi_i}{\partial y}\\[3mm]
\frac{\partial \varphi_i}{\partial z}\end{array}\right] = \bfJ^{-1}(\xi,\eta,0) \left[\begin{array}{>{\displaystyle}c}
\frac{\partial \varphi_i}{\partial \xi}\\[3mm]
\frac{\partial \varphi_i}{\partial \eta}\\[3mm]
0\end{array}\right] ,\quad \bfJ(\xi,\eta,0) := \left[\begin{array}{>{\displaystyle}c>{\displaystyle}c>{\displaystyle}c}
\frac{\partial x^h}{\partial \xi} & \frac{\partial y^h}{\partial \xi} & \frac{\partial z^h}{\partial \xi}\\[3mm]
\frac{\partial x^h}{\partial \eta} & \frac{\partial y^h}{\partial \eta} & \frac{\partial z^h}{\partial \eta}\\[3mm]
n_{Tx} & n_{Ty} & n_{Tz}\end{array}\right] .
\]
This low order approximation has the advantage of yielding a constant Jacobian from a linear approximation. For some applications this
is, however, 
offset
by the problem of having a discontinuous normal between elements.
\end{rem}

\subsection{Finite element formulation}

We can now introduce finite element spaces constructed from the basis previously discussed by defining 
\begin{equation}\label{spacevA}
W^h := \{ v: {v\vert_T \circ\bfF\in P^k(\hat T),\; \forall T\in\mathcal{T}_h};\; v\; \in C^0(\Sigma_h)\},
\end{equation}
(in the numerics, we use $k=1$), and the finite element method reads: Find $\bfu_h\in V^h := [W^h]^3$ such that 
\begin{equation}
a_{\Sigma_h}(\bfu_h,\bfv) = l_{\Sigma_h}(\bfv) ,\quad \forall \bfv \in V^h,
\end{equation}
where 
\begin{align}\nonumber
a_{\Sigma_h}(\bfu,\bfv)   = {} & (2\mu\bfeps_{\Sigma_h}(\bfu), \bfeps_{\Sigma_h}(\bfv))_{\Sigma_h}-(4\mu\bfeps_{\Sigma_h}(\bfu)\cdot\bfn^h, \bfeps_{\Sigma_h}(\bfv)\cdot\bfn^h)_{\Sigma_h}\\\nonumber
& + ({\lambda_0}\, \nabla_{\Sigma_h}\cdot\bfu,\nabla_{\Sigma_h}\cdot\bfv)_{\Sigma_h}\\\nonumber
\end{align}
and
$l_{\Sigma_h}(\bfv) = (\bff,\bfv)_{\Sigma_h}$. 

\subsection{Extension to surfaces with a boundary}

If the surface $\Sigma$ has a 
boundary $\partial \Sigma$ we assume that $\partial \Sigma  = \cup_i \partial\Sigma_i$ where $\partial\Sigma_i$ 
are closed components. On each of the components $\partial\Sigma_i$ we let $\bfq_j : \partial\Sigma_i \rightarrow {\bf R}^3, j=1,2,3,$ be smooth orthonormal vector fields. We strongly impose  homogeneous Dirichlet 
boundary conditions of the type
\begin{equation}
\bfq_j \cdot \bfu = 0 \quad \text{on $\partial \Sigma_i$}, \quad 1 \leq j \leq d_i,
\end{equation}
where $d_i = 1,2,$ or $3$, and weakly the remaining Neumann condition 
\begin{equation}
( \bfn_{\partial \Sigma} \cdot \bfsig_\Sigma^P (\bfu ) ) \cdot \bfq_j = 0  , \quad d_i < j \leq 3,
\end{equation}
where $\bfn_{\partial \Sigma}$ is the unit vector that is normal to $\partial \Sigma$ and tangent to $\Sigma$. Note that not every combination of boundary conditions and right hand side leads to a well 
posed problem.


\section{Numerical examples}

In the numerical examples below, the geometry is represented by flat facets, and the normals are taken as the exact normal in the nodes, interpolated linearly inside each element. Our experience is that similar results are obtained if we use $L_2-$projections of the 
flat facet normals in the nodes and then interpolate these linearly.

\subsection{Pulling a cylinder}
\begin{figure}[h]
\begin{center}
\includegraphics[width=5in]{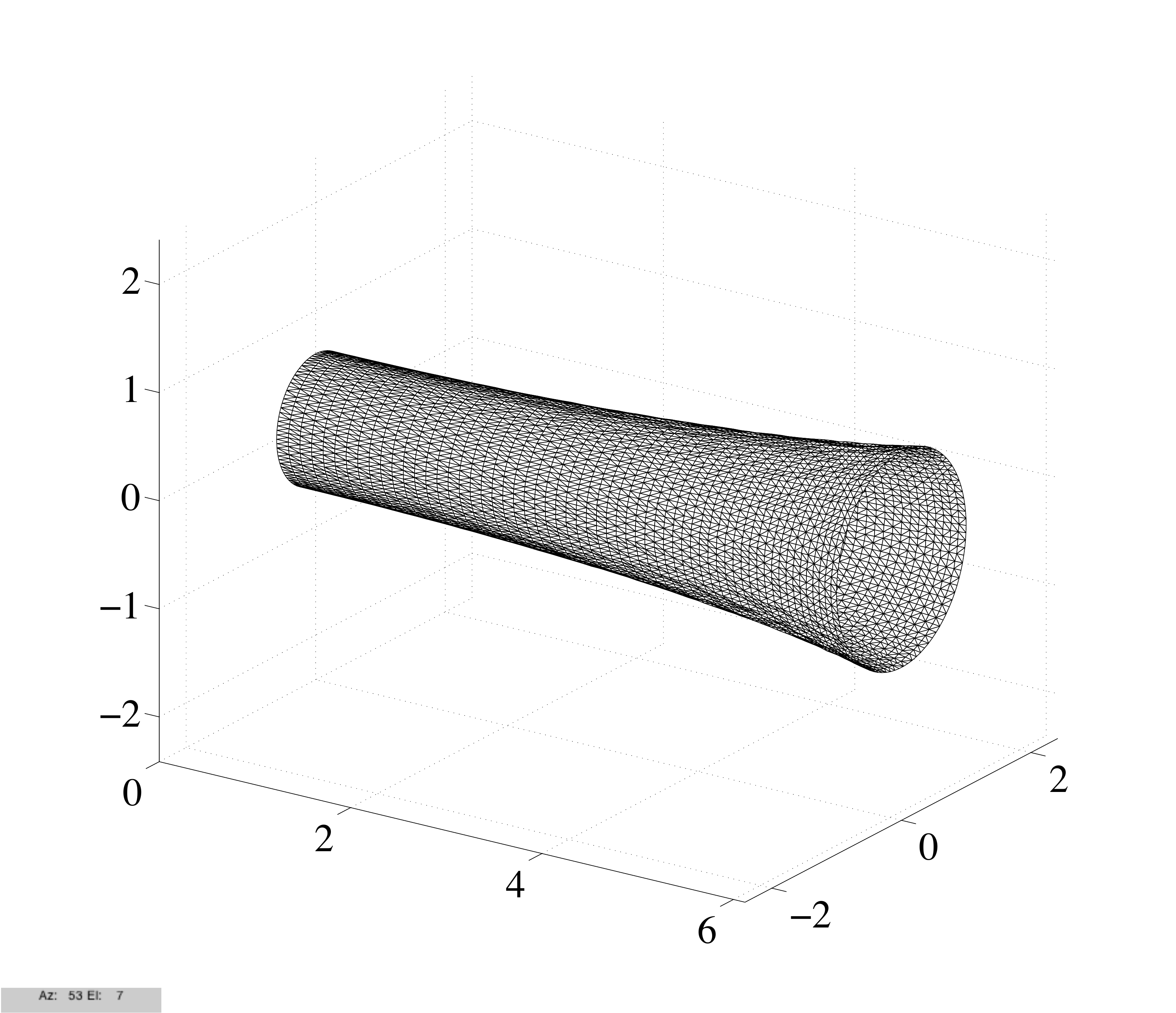}
\end{center}
\caption{Displacements (exaggerated by one order of magnitude) on a particular mesh.}\label{fig:unstab}
\end{figure}

\begin{figure}[h]
\begin{center}
\includegraphics[width=5in]{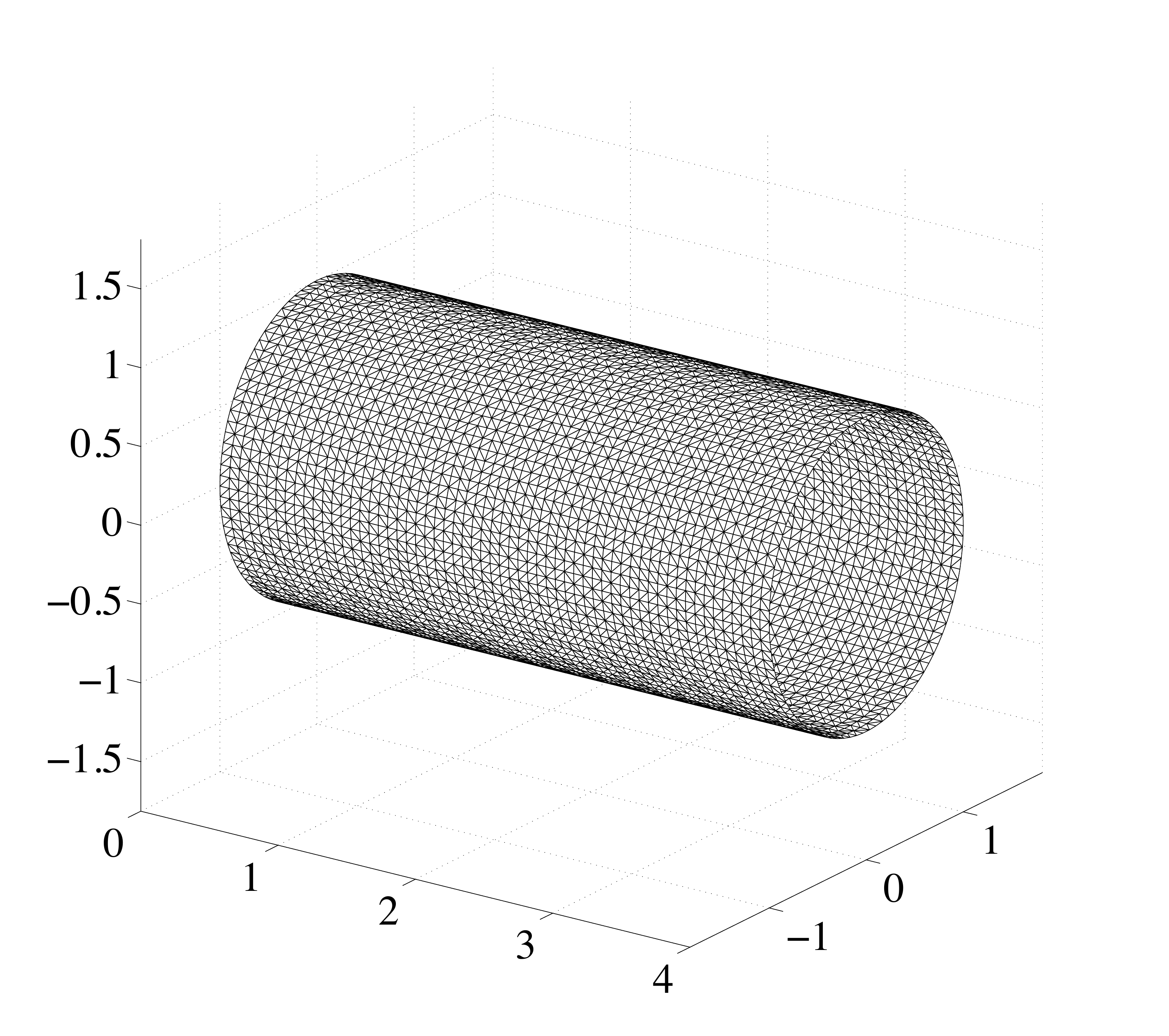}
\end{center}
\caption{The cylinder before deformation.}\label{fig:cal}
\end{figure}

\begin{figure}[h]
\begin{center}
\includegraphics[width=5in]{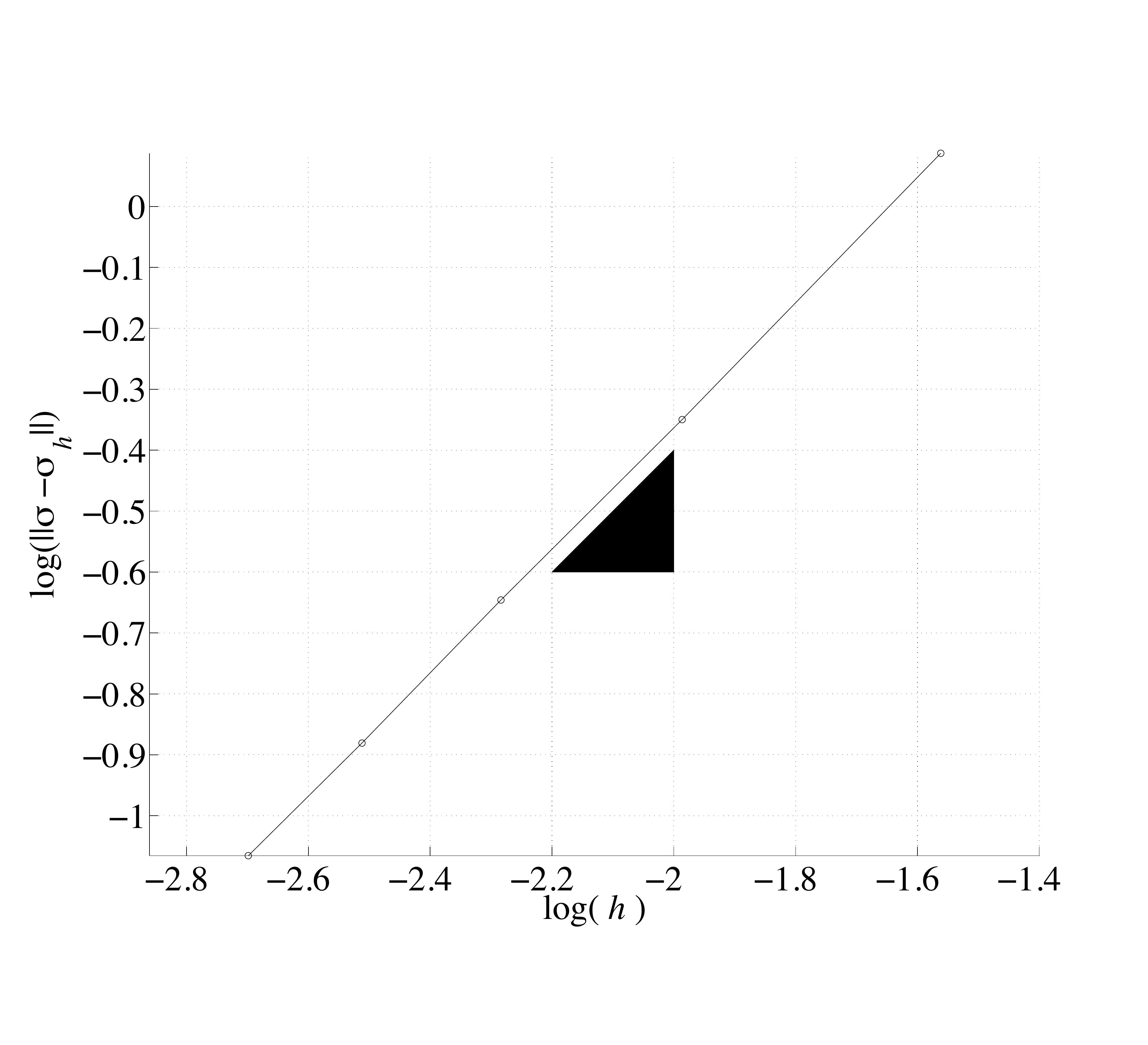}
\end{center}
\caption{Stress convergence for the cylinder}\label{fig:stresserr}
\end{figure}

We consider a cylindrical shell of radius $r$ and thickness $t$, with open ends at $x=0$ and at $x=L$, and with fixed longitudinal displacements at $x=0$, and radial at $x=L$, carrying a horizontal surface load per unit area
\[
f(x,y,z) = \frac{F}{2\pi r}\frac{x}{L^2} ,
\]
where $F$ has the unit of force. The resulting longitudinal stress
is 
\[
\sigma = \frac{F\left(\displaystyle 1-(x/L)^2\right)}{4\pi r t}.
\]
We take as an example a cylinder of radius $r=1$ and length $L=4$, with material data $E=100$ and $\nu = 1/2$, with thickness $t=10^{-2}$, and with $F=1$.
In Fig. \ref{fig:unstab} we show the solution (exaggerated 10 times) on a particular mesh (shown in Fig. \ref{fig:cal}).
Note that the lateral contraction creates radial displacements depending on the size of stress. Finally, in Fig. \ref{fig:stresserr} we show the $L_2$ error in stresses, $\Vert \bfsig -\bfsig_h\Vert_{L_2(\Omega)}$, where $\bfsig := \bfsig^P_\Sigma(\bfu)$ and $\bfsig_h := \bfsig^P_\Sigma(\bfu_h)$, which shows the expected first order convergence for our $P^1$ approximation. The black triangle shows the 1:1 slope.

\subsection{A torus with internal pressure}

\begin{figure}[h]
\begin{center}
\includegraphics[width=5in]{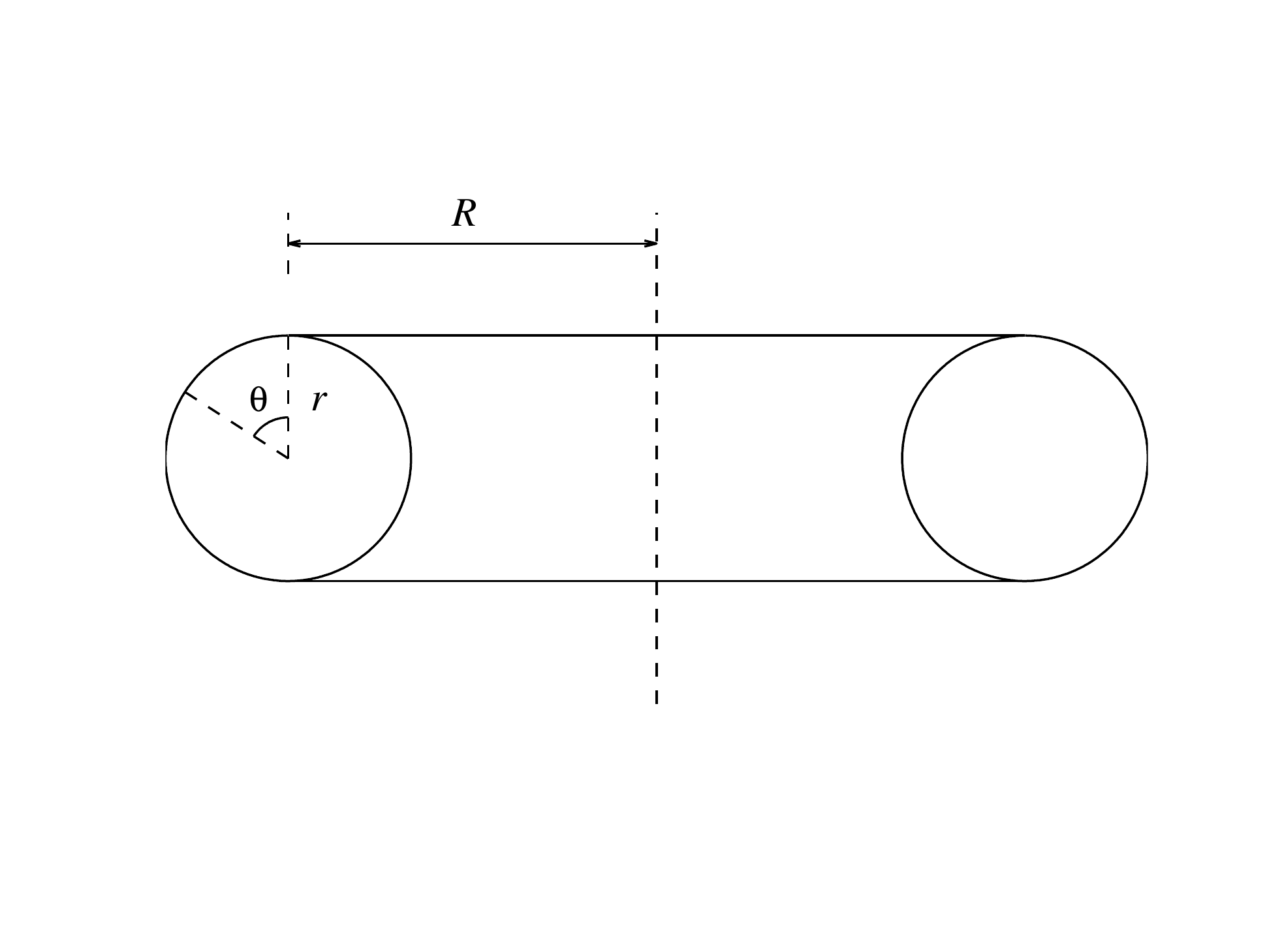}
\end{center}
\caption{A cut through of the torus}\label{fig:torus}
\end{figure}

We consider a torus with internal gauge pressure $p$ for which the stresses are statically determinate. 
Using the angle and radii defined in Fig. \ref{fig:torus}, the principal stresses are given by
\[
\sigma_1 =\frac{pr}{2t}, \; \sigma_2 = \frac{pr}{t}\left(1-\frac{r\sin{\theta}}{2(R+r\sin{\theta})}\right) ,
\]
where $\sigma_1$ is the longitudinal stress, $\sigma_2$ the hoop stress, and $t$ is the thickness of the surface of the torus. The constitutive parameters and thickness where chosen as in the cylinder example, and
we set $R=1$, $r=1/2$, and $p=1$.

Again we compute the stress error $\| \bfsig-\bfsig_h\|_{L_2(\Omega)}$. We show the observed convergence in Fig. \ref{fig:toruserr} at a rate of about $3/4$ (the slope of the black triangle), which is suboptimal, but does occur in
problems where elliptic regularity is an issue, cf. \cite{buha}, Lemma 10. We thus attribute 
this loss of convergence to the
load now being in the normal direction of the shell, for which we do not have ellipticity.

\newpage
\begin{figure}[h]
\begin{center}
\includegraphics[width=5in]{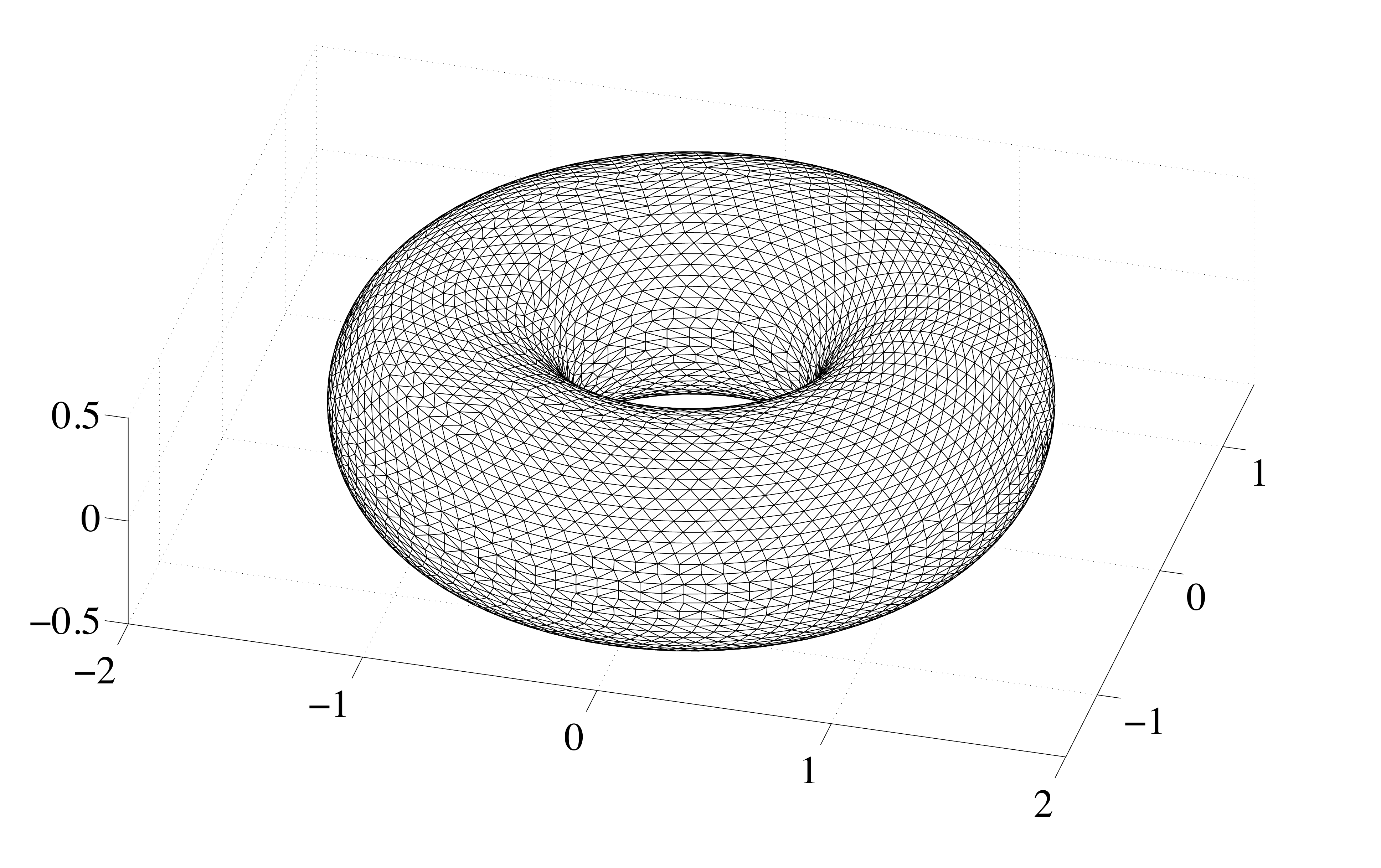}
\end{center}
\caption{A typical mesh on the torus.}\label{fig:toruserr}
\end{figure}
\begin{figure}[h]
\begin{center}
\includegraphics[width=5in]{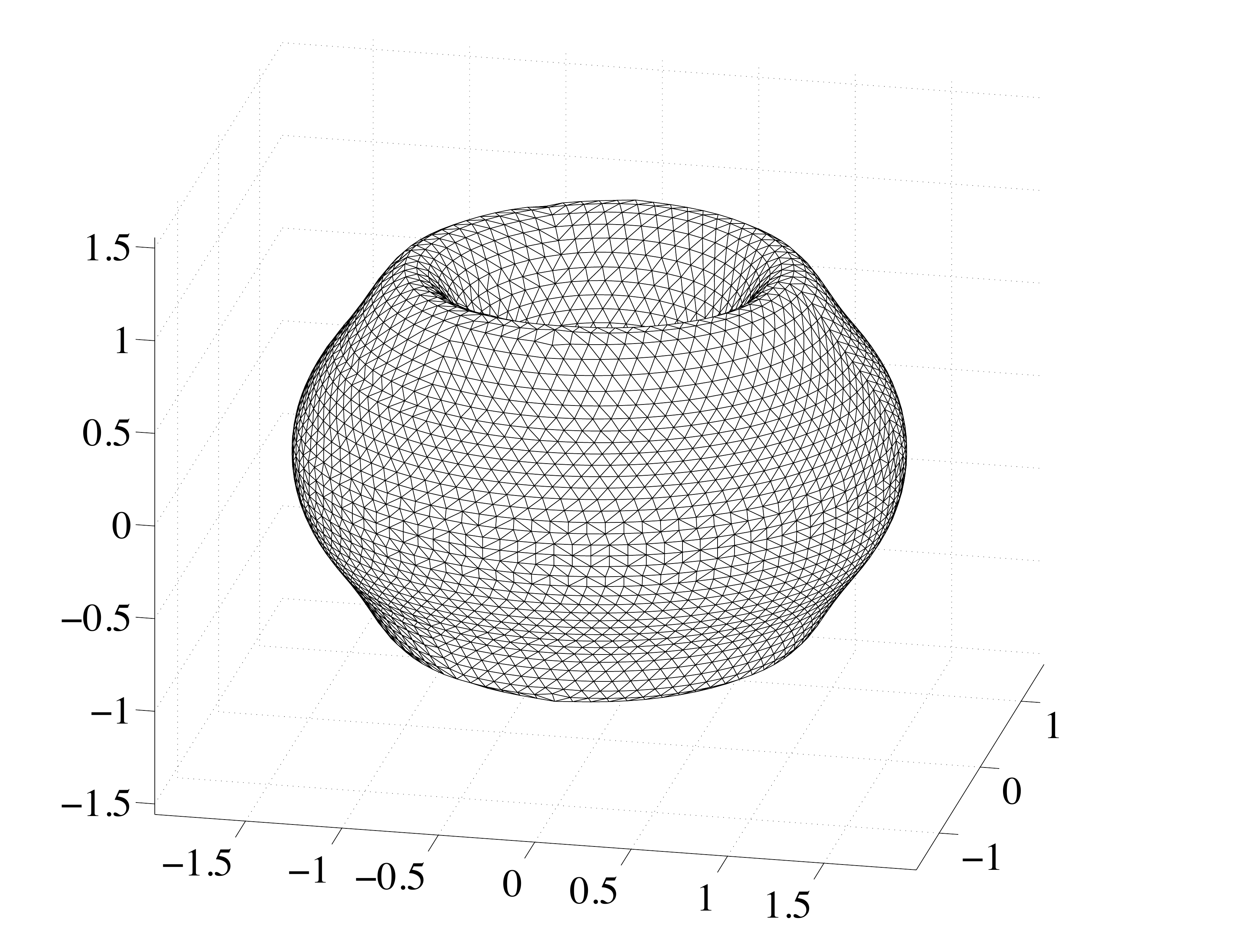}
\end{center}
\caption{Deformations in the torus case, exaggerated by two orders of magnitude.}\label{fig:toruserr}
\end{figure}
\begin{figure}[h]
\begin{center}
\includegraphics[width=5in]{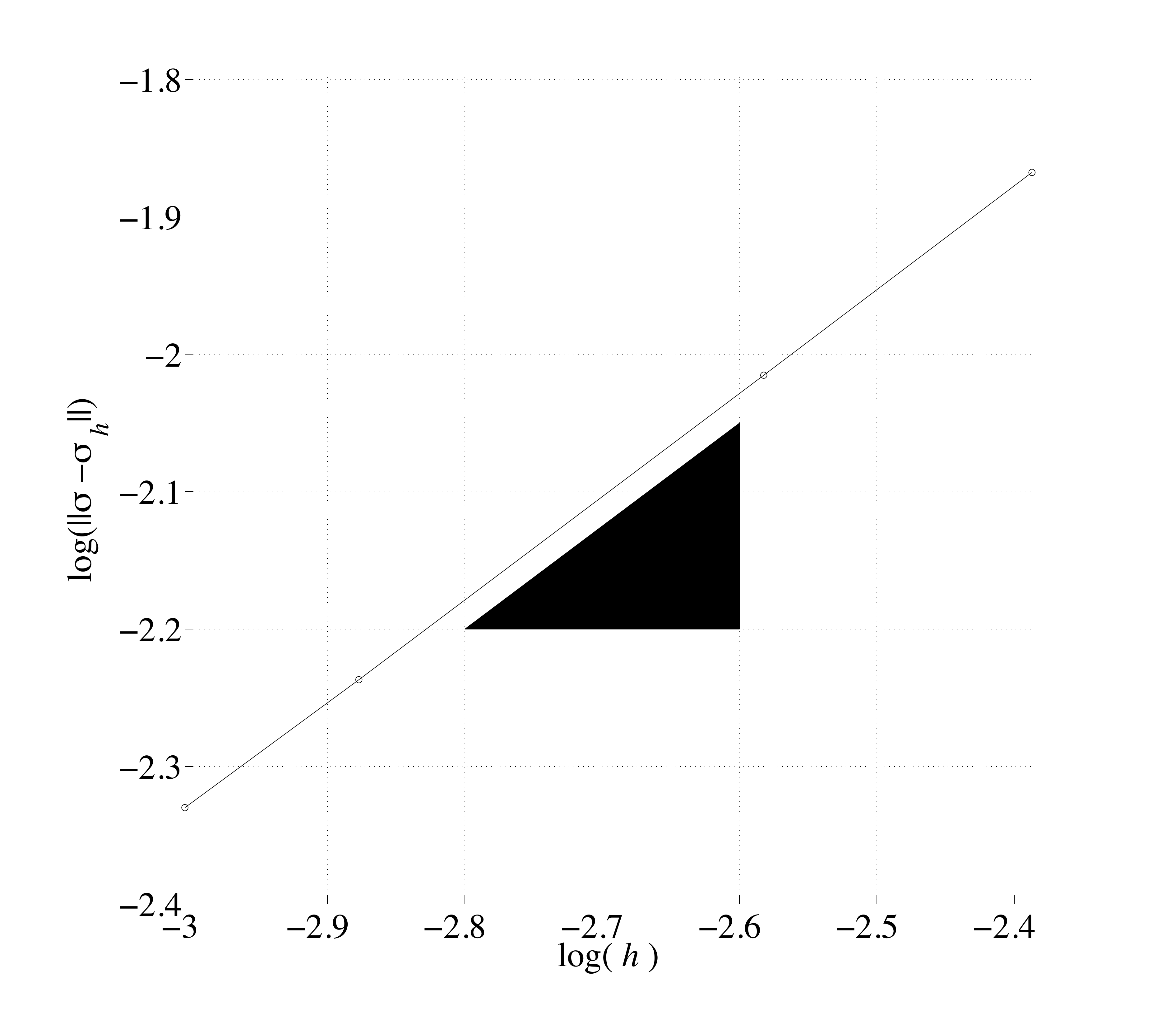}
\end{center}
\caption{Convergence of the stresses in the torus case.}\label{fig:toruserr}
\end{figure}

\end{document}